\documentclass[11pt,a4paper]{article}

\usepackage{epsf,epsfig,amsfonts,amsgen,amsmath,amstext,amsbsy,amsopn,amsthm
}
\usepackage{amsmath,times,mathptmx}
\usepackage{amsfonts,amsthm,amssymb}
\usepackage{amsfonts}
\usepackage{graphics}
\usepackage{latexsym,bm}
\usepackage{amsfonts,amsthm,amssymb,bbding}
\usepackage{indentfirst}
\usepackage{graphicx}
\usepackage{color}
\usepackage[colorlinks=true,anchorcolor=blue,filecolor=blue,linkcolor=blue,urlcolor=blue,citecolor=blue]{hyperref}
\usepackage{float}
\usepackage{tikz}
\newtheorem{thm}{Theorem}

\newtheorem{lemma}{Lemma}
\newtheorem{false statement}{False statement}

\theoremstyle{definition}

\newtheorem{claim}{Claim}

\newtheorem{case}{Case}
\newtheorem{subcase}{Subcase}[case]

\begin{document}

\title{Extremal distance spectra of graphs and essential connectivity
\footnote{Supported by Natural Science Foundation of Xinjiang Uygur Autonomous Region (No. 2024D01C41), NSFC (Nos. 12361071 and 11901498).}}
\author{{Daoxia Zhang, Dan Li}\thanks{Corresponding author. E-mail: ldxjedu@163.com.}, {Wenxiu Ding}\\
{\footnotesize College of Mathematics and System Science, Xinjiang University, Urumqi 830046, China}}
\date{}

\maketitle {\flushleft\large\bf Abstract:}
A graph is non-trivial if it contains at least one nonloop edge. The essential connectivity of $G$, denoted by  $\kappa'(G)$, is the minimum number of vertices of $G$ whose removal produces a disconnected graph with at least two components are non-trivial. In this paper, we determine  the $n$-vertex graph of given essential connectivity with minimum distance spectral radius. We also characterize the extremal graphs attaining the minimum distance spectral radius among all connected graphs with fixed essential connectivity and minimum degree. Furthermore, we characterize the extremal digraphs with minimum distance spectral radius among the strongly connected digraphs with given essential connectivity.
\vspace{0.1cm}
\begin{flushleft}
\end{flushleft}

\section{Introduction}
All graphs considered in this paper are simple and connected. Let $G=(V(G),E(G))$ be a graph with vertex set $V(G)$ and edge set $E(G)$. The distance between $v_{i}$ and $v_{j}$, denote by $d_{ij}(G)$, is the length of a shortest path from $v_{i}$ to $v_{j}$. The distance matrix of $G$, denote by $D(G)= (d_{i j})$, is the $|V(G)|$ × $|V(G)|$ matrix whose $(i, j)$-entry is $d_{ij}$. Clearly, $D(G)$ is a real symmetric matrix with zeros on the diagonal. Let $\lambda_i(D(G))$ denote the $i$th largest eigenvalue of $D(G)$. Particularly, the largest eigenvalue of $D(G)$ is called the distance spectral radius of $D(G)$. By the Perron-Frobenius theorem, $\lambda _{1}(D(G))$ is always positive and there exists a unique positive unit eigenvector $X=(x_{1},x_{2}, \ldots ,x_{n})^{T}$ corresponding to $\lambda _{1}(D(G))$ such that $D(G)X=\lambda_{1}(D(G))X$, which is called the Perron-vector of $D(G)$.   
A graph can be used to represent the topology of an interconnected network, where the vertices in the graph correspond to the processors of the network, and the edges correspond to the connections between the communication processors. As we all know, the connectivity is an important parameter for measuring the network reliability. The greater the connectivity, the greater the reliability of the network. We use $\kappa(G)$ to represent the connectivity of a graph. In 1983, Harary \cite{F. Harary} introduced the concept of conditional connectivity by imposing some conditions on the components of $G-S$, where $S$ is a subset of edges or vertices. A graph is non-trivial if it contains at least one nonloop edge. A vertex-cut $S$ of $G$ is essential if at least two components of $G-S$ are non-trivial. In particular, the essential connectivity of $G$, denote by $\kappa'(G)$ or $\kappa'$, is the size of a minimum essential vertex-cut. As a conditional connectivity, essential connectivity plays a significant role in the robustness analysis of computer networks.
 
Over the past several decades, the distance spectral radius of graphs with some given parameters have been well investigated. In 2010, Liu \cite{Z.Z. Liu} characterized graphs with minimum distance spectral radius among $n$-vertex graphs with fixed vertex connectivity (resp. matching number, chromatic number). Zhang \cite{X.L. Zhang2} determined the $n$-vertex graphs of given diameter with the minimum distance spectral radius. Subsequently, Zhang, Li and Gutman \cite{M.J. Zhang} generalized this result of \cite{X.L. Zhang2} by characterizing the graphs of order $n$ with given connectivity and diameter having minimum distance spectral radius. Moreover, they determined the minimum distance spectral radius of graphs among the $n$-vertex graphs with given connectivity and independence number. In 2021, Zhang and Lin \cite{Y.K. Zhang} presented sufficient conditions in terms of the distance spectral radius to guarantee the existence of a perfect matching in graphs and bipartite graphs, respectively, and then Zhang, Lin, Liu and Zheng \cite{Y.K. Zhang2} generalized the result of \cite{Y.K. Zhang}. In addition, they determined the extremal graph attaining the minimum distance spectral radius among all bipartite graphs with a unique perfect matching. 
Fan, Lin and Lu \cite{D.D. Fan} provided sufficient conditions in terms of the spectral radius for a graph to be 1-tough with minimum degree $\delta$ and $t$-tough with $t\geq1$ being an integer, respectively. Recently, Lou, Liu and Shu \cite{J. Lou} gave sufficient conditions based on the distance spectral radius to guarantee a graph to be 1-tough with minimum degree $\delta$. In addition, they studied sufficient conditions with respect to the distance spectral radius for a graph to be $t$-tough, where $t$ or $\frac{1}{t}$ is a positive integer. For more results about the relationships between structural properties and the distance spectral radius of graphs, we refer the readers to {\cite{S.S. Bose, S.S. Bose1, Y.Y. Chen, X.L. Zhang1}}.

Motivated by the aforementioned results of the distance spectral extremum problems for given graph parameters, we are interested in investigating distance spectral extrema with regard to essential connectivity of graphs and our result is as follows. As usual, $K_{n}$ denotes the complete graph of order $n$. For two vertex-disjoint graphs $G_{1}$ and $G_{2}$, we denote by $G_{1}\cup G_{2}$ the disjoint union of $G_{1}$ and $G_{2}$. The join $G_{1}\vee G_{2}$ is the graph obtained from $G_{1}\cup G_{2}$ by adding all edges between $V(G_{1})$ and $V(G_{2})$.

\vspace*{2mm}
\begin{thm}\label{thm1.1}
Let $G$ be a connected graph of order $n\geq\kappa'+4$ with essential connectivity $\kappa'$. Then $$\lambda_{1}(D(G))\geq\lambda_{1}(D(K_{\kappa'}\vee(K_{2}\cup K_{n-\kappa'-2}))),$$ 
with equality if and only if $G\cong K_{\kappa'}\vee(K_{2}\cup K_{n-\kappa'-2})$.
\end{thm}

Let $\mathcal{G}_{n}^{\kappa', \delta}$ be the set of graphs with given order $n$, minimum degree $\delta$ and essential connectivity $\kappa'$. Notice that there is no bound between $\kappa'$ and $\delta$. If $\kappa'>\delta-1$, let $G_{n}^{\kappa', \delta}$ be the graph obtained from $\{z\}\cup(K_{\kappa'}\vee (K_{1}\cup K_{n-\kappa'-2}))$ by adding $\delta-1$ edges between the isolated vertex $z$ and $K_{\kappa'}$ and adding an edge between $z$ and $K_{1}$. If $\kappa'\leq\delta-1$, let $G_{n}^{\kappa',\delta}= K_{\kappa'}\vee (K_{n-\delta-1}\cup K_{\delta-\kappa'+1})$. Based on Theorem \ref{thm1.1}, we also study the problem of minimum distance spectral radius of graphs with given essential connectivity and minimum degree. Then we have the following result.

\vspace*{2mm}
\begin{thm}\label{thm1.2}
Let $G$ be a connected graph of order $n \geq \kappa' + 4$ with minimum degree $\delta$ and essential connectivity $\kappa'$. Then 
$$\lambda_{1}(D(G))\geq\lambda_{1}(D(G_{n} ^{\kappa', \delta})),$$
with equality if and only if $G\cong G_{n} ^{\kappa', \delta}.$
\end{thm}

Let $D(\overrightarrow{G})=(d_{ij})$ be the distance matrix of a digraph $\overrightarrow{G}$, where $d_{ij}$ is the length of a shortest path from $v_{i}$ and $v_{j}$. A digraph $\overrightarrow{G}$ is strongly connected if for every pair $u,v\in V(\overrightarrow{G})$, there exists a directed path from $u$ to $v$ and a directed path from $v$ to $u$. The matrix $D(\overrightarrow{G})$ is nonnegative and irreducible when $\overrightarrow{G}$ is strongly connected. The eigenvector of $D(\overrightarrow{G})$ with the largest modulus is called the distance spectral radius of the digraph $\overrightarrow{G}$, denote by $\lambda_{1}(D(\overrightarrow{G}))$. For more results on the distance matrix and its spectral properties, we refer the readers to the survey \cite{M. Aouchiche}. The essential connectivity of a strongly connected digraph is the minimum number of vertices whose removal produces a digraph with at least two components are strongly connected and non-trivial. 

We also consider an analogue for digraphs. In 2012, Lin, Yang, Zhang and Shu \cite{H.Q. Lin} characterized the extremal digraphs with minimum distance spectral radius among all digraphs with given vertex connectivity and the extremal graphs with minimum distance spectral radius among all graphs with given edge connectivity. Moreover, they gave the exact value of the distance spectral radius of those extremal digraphs and graphs. They also characterized the graphs with the maximum distance spectral radius among all graphs of fixed order with given vertex connectivity 1 and 2. Lin and Shu \cite{H.Q. Lin2} first gave sharp upper and lower bounds for the distance spectral radius for strongly connected digraph, then they characterized the digraph having the maximal and minimal distance spectral radius among all strongly connected digraph. In addition, they determined the extremal digraph with the minimal distance spectral radius with given arc connectivity and dichromatic number. In 2021, Xi, So and Wang \cite{W.G. Xi} determined the strongly connected digraphs minimizing $\lambda_{1}(D(\overrightarrow{G}))$ among all strongly connected digraphs with fixed order $n$ and diameter $d$ for $d=1,2,3,4,5,6,7,n-1$.

Let $\overrightarrow{\mathcal{D}}_{n,k}$ be the set of strongly connected digraphs with given essential connectivity $\kappa'(\overrightarrow{G})=k$. Denote by $\overrightarrow{G}_{1}\nabla \overrightarrow{G}_{2}$ the digraph obtained from two disjoint digraphs $\overrightarrow{G}_{1}$ and $\overrightarrow{G}_{2}$ with vertex set $V(\overrightarrow{G}_{1})\cup V(\overrightarrow{G}_{2})$ and arc set $E=E(\overrightarrow{G}_{1})\cup E(\overrightarrow{G}_{2})\cup\{uv, vu \mid u\in V(\overrightarrow{G}_{1}), v\in V(\overrightarrow{G}_{2})\}$. 
The complete digraph of order $n$ is the digraph $\overrightarrow{K_{n}}$ in which every pair of vertices is connected by a pair of opposite arcs.
Let $\overrightarrow{G}_{n}^{k,n_{1}}=\overrightarrow{K}_{k}\nabla (\overrightarrow{K}_{n_{1}}\cup \overrightarrow{K}_{n-k-n_{1}})\cup F$, where $F=\{uv \mid u\in V(\overrightarrow{K}_{n_{1}}), v\in V(\overrightarrow{K}_{n-k-n_{1}})\}$. Let $\overrightarrow{\mathcal{G}}(n,k)=\{\overrightarrow{G}_{n}^{k,n_{1}} \mid 2 \leq n_{1}\leq n-k-2\}$. Clearly, $\overrightarrow{\mathcal{G}}(n,k)\subseteq \overrightarrow{\mathcal{D}}_{n,k}$. Thus we have the following results.

\vspace*{2mm}
\begin{thm}\label{thm1.3}
The digraphs $\overrightarrow{G}_{n}^{k,2}$ and $\overrightarrow{G}_{n}^{k,n-k-2}$ minimize the distance spectral radius among all digraphs in $\overrightarrow{\mathcal{D}}_{n,k}$.
\end{thm}

\begin{figure}[H]
\centering
\includegraphics[width=0.7\linewidth]{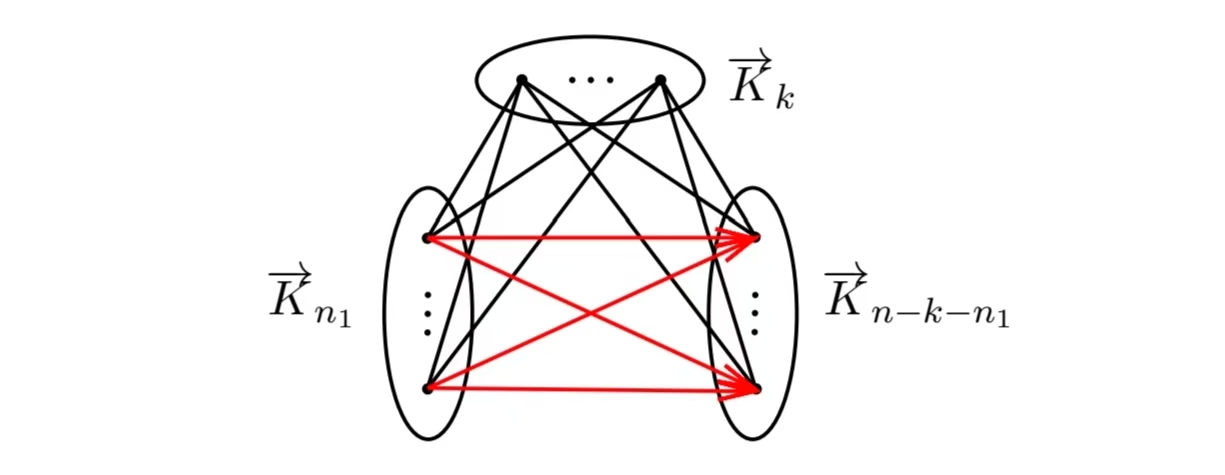}
\caption{The digraph $\overrightarrow{G}_{n}^{k,n_{1}}$.}
\end{figure}

This article is organized as follows: In Section \ref{2}, we characterize the graph with minimum distance spectral radius among all the connected graphs with given $n$ vertices and essential connectivity. Based on this, we determine the $n$-vertex graphs of fixed essential connectivity and minimum degree with the minimum distance spectral radius among all connected graphs. In Section \ref{3}, we characterize the extremal digraphs with the minimum distance spectral radius among the strongly connected digraphs with given essential connectivity. 

\section{Proof of Theorem \ref{thm1.1} and \ref{thm1.2}}\label{2}

\begin{lemma}\cite{C.D. Godsil}\label{lem1.1}
Let $e$ be an edge of $G$ such that $G-e$ is connected. Then $\lambda_{1}(D(G)) < \lambda_{1}(D(G-e))$.
\end{lemma}

\begin{lemma}\cite{J. Lou}\label{lem1.2}
Let $n$, $c$, $s$, $p$ and $n_{i}$ ($1\leq i\leq c$) be positive integers with $n_{1}\geq 2p$, $n_{1}\geq n_{2}\geq \cdots \geq n_{c}\geq p$ and $n_{1}+n_{2}+\cdots +n_{c}=n-s$. Then $$\lambda_{1}(D(K_{s}\vee(K_{n_{1}}\cup K_{n_{2}}\cup \cdots \cup K_{n_{c}}))\ge \lambda_{1}(D(K_{s}\vee(K_{n-s-p(c-1)}\cup(c-1)K_{p}))),$$
with equality if and only if $(n_{1}, n_{2}, \ldots, n_{c})=(n-s-p(c-1), p, \ldots, p)$.
\end{lemma}

\noindent{\bf{Proof of Theorem \ref{thm1.1}}.}
We assume that $G$ is the connected graph with minimum distance spectral radius among all graphs with given order $n$ and essential connectivity $\kappa'$. Let $S$ be the subset of $V(G)$ with $|S|=\kappa'$ such that $G-S$ contains at least two non-trivial components and $A_{1}, A_{2}, \ldots, A_{l}$ be the $l$ components of $G-S$, where $|V(A_{i})|=n_{i}$ for $i=1, 2, \ldots, l$ and $l\geq2$. By the minimality of $\lambda_{1}(D(G))$, we can obtain that $G\cong K_{\kappa'}\vee (K_{n_{1}}\cup K_{n_{2}}\cup \cdots \cup K_{n_{l}})$. We first assert that $l=2$. Otherwise, $l\geq3$. Without loss of generality, we suppose that $A_1$ is one of the non-trivial components. Let $G'=K_{\kappa'}\vee (K_{n_{1}}\cup K_{n_{2}'})$, where $n_{2}'=\sum_{i=2}^{l}n_{i}$. It is easy to see that $\lambda_{1}(D(G'))<\lambda_{1}(D(G))$ by Lemma \ref{lem1.1}, a contradiction. This means that $l=2$. Hence, $G-S$ has exactly two non-trivial components $K_{n_{1}}$ and $K_{n_{2}}$. Then $G\cong K_{\kappa'}\vee(K_{n_{1}}\cup K_{n_{2}})$, where $\text{min}\{n_1,n_2\}\geq 2$. Without loss of generality, we assume that $n_{1}\leq n_{2}$. Combining this with Lemma \ref{lem1.2}, we have
$$\lambda_{1}(D(K_{\kappa'}\vee (K_{2}\cup K_{n-\kappa'-2})))\leq \lambda_{1}(D(K_{\kappa'}\vee(K_{n_{1}}\cup K_{n_{2}}))),$$
with equality if and only if $(n_{1},n_{2})=(2,n-\kappa'-2)$. By the minimality of $\lambda_{1}(D(G))$, we deduce that $G\cong K_{\kappa'}\vee (K_{2}\cup K_{n-\kappa'-2})$.

This completes the proof. 
\begin{flushright}
$\square$
\end{flushright}

Let $N_{G}(u)=\left \{ v\in V(G)\mid uv\in E(G) \right \}$ be the set of neighbors of $u$ in $G$ and $N_{G}[u]=N_G(u)\cup\{u\}$. The degree of the vertex $u$ is $d_{G}(u)=|N_{G}(u)|$. In this paper, we always use $J$ to denote the all-one matrix and $O$ to denote the zero matrix. Now, we are ready to present the proof of Theorem \ref{thm1.2}.

\noindent{\bf{Proof of Theorem \ref{thm1.2}}.}
Suppose that $G$ is a connected graph attaining the minimum distance spectral radius in $\mathcal{G}_{n}^{\kappa',\delta}$. By the definition of essential connectivity, there exists some nonempty subset $S\subseteq  V(G) $ with $|S| =\kappa'$ such that $G-S$ contains at least two non-trivial components. Let $A_{1},A_{2}, \dots , A_{l}$ be  the components of $G-S$ and at least two of them be non-trivial, and $|V(A_{i})| =n_{i}$ for $1\le i\le l$. Based on the relationship between $\kappa'$ and $\delta$, we will divide the proof into the following two cases.

\begin{case}\label{case1.1}   
$\kappa '> \delta -1$.
\end{case}

For any vertex $v\in V(G)$ and any subset $S\subseteq  V(G) $, let $N_{S}(v)=N_{G}(v)\cap S$ and $d_{S}(v) =\left | N_{G}(v)\cap S \right |$. Choose a vertex $z\in V(G)$ such that $d_{G}(z) =\delta$ and $d_{S}(z)=b$.

\begin{claim}\label{claim1.1}
$z\notin S$.     
\end{claim}

Otherwise, $z\in S$. Firstly, we can obtain that there are at least two $i$ of $\{1, 2, \ldots, l\}$ such that $N_{G}(z)\cap V( A_{i})\neq \emptyset$. Otherwise, if $N_{G}(z)\cap V (A_{i})= \emptyset$ for any $1\le i\le l$, then we can get a new subset $S_{1}=S-z $ with $|S_{1}|=\kappa '-1$ such that $G-S_{1}$ contains $l+1$ components $z, A_{1}, A_{2}, \ldots, A_{l}$, where at least two of them  are non-trivial components, contradicting the definition of essential connectivity $\kappa'$ of $G$. If there is only one positive integer $i\in [1,l]$ such that $N_{G}(z) \cap V(A_{i} )\neq\emptyset$, then we can also obtain a new subset $S_{2}=S-z$ with $|S_{2}|=\kappa'-1$ such that $G-S_{2}$ contains $l$ components $A_{1}, \ldots, A_{i}\cup\{z\}, \ldots, A_{l}$, where at least two of them are non-trivial components, this still leads to a contradiction. Let $P_{i} =N_{G} (z)\cap V(A_{i} )$ and $|P_{i}| =p_{i}$ for $1\le i\le l$. Then $|V(A_{i})\setminus P_{i}|=n_{i}-p_{i}$ and $\sum_{i=1}^l p_i=\delta-b$. Without loss of generality, we suppose that $A_{1}$ is one of the non-trivial components. By the minimality of $\lambda_{1}(D(G))$, we can obtain that $G-z \cong K_{\kappa '-1 } \vee (K_{n_{1}}\cup K_{n_{2}}\cup \cdots \cup K_{n_{l}})$. It is clear that $l=2$. Otherwise, we can construct a new graph $H$ obtained from $G$ by joining $K_{n_{2}}, K_{n_{3}},  \ldots , K_{n_{l}}$. Obviously, $H\in \mathcal{G}_{n}^{\kappa', \delta}$ and $\lambda_{1}( D(H)) < \lambda_{1}(D(G))$ by Lemma \ref{lem1.1}, a contradiction. Thus, $G-z\cong K_{\kappa '-1}\vee(K_{n_{1}} \cup K_{n_{2}})$. Based on the relationship between $p_{i}$ and $n_{i}$ for $i=1, 2$, we will discuss the case into two subcases.
 
\begin{subcase}\label{subcase1.1}
$p_1<n_1$ or $p_2<n_2$.     
\end{subcase}

Without loss of generality, assume that $p_1<n_1$. Then $1\leq p_1<n_1$ and $1\leq p_2\leq n_2$. Next, we can construct a new graph $G_{1}$ from $G$ through the following steps:
 
\text{\bf{Step 1:}} Remove $p_{2}$ edges between $z$ and $K_{p_{2}} $;
 
\text{\bf{Step 2:}} Add $p_{2}$ edges between $z$ and $S\setminus N_{S}[z]$;

\text{\bf{Step 3:}} Choose a vertex $w$ in $V(K_{n_{1}-p_{1}})$ and join $w$ to $V(K_{n_{2}})$.
\begin{figure}[H]
\centering
\includegraphics[width=0.8\linewidth]{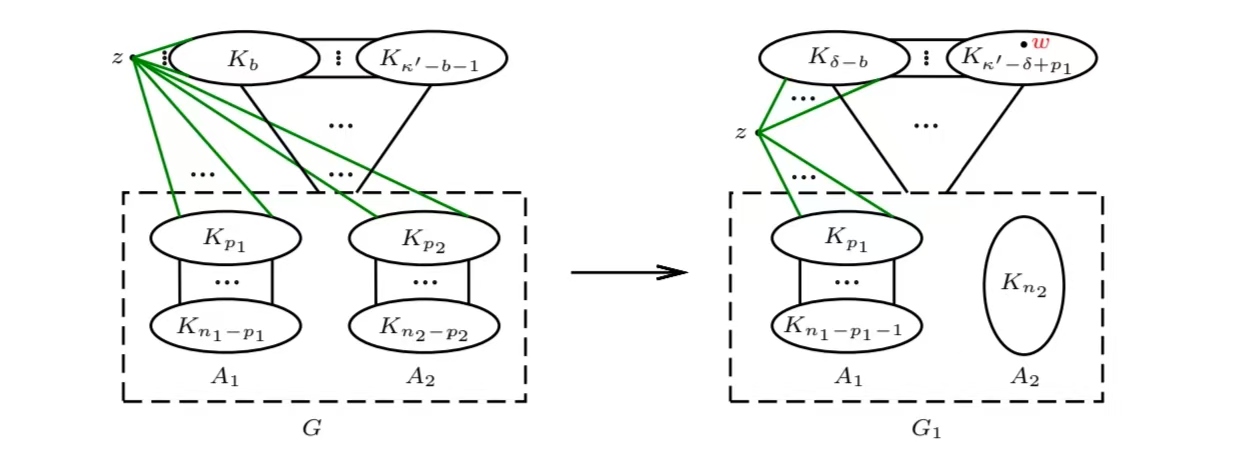}
\caption{The graphs $G$ and $G_{1}$.}
\end{figure}

Then $G_{1}-z\cong K_{\kappa '}\vee (K_{n_{1}-1}\cup K_{n_{2}})$ and $G_1\in \mathcal{G}_{n}^{\kappa', \delta}$.
Let $X$ be the Perron vector of $D(G_{1})$ and $x(v)$ denote the entry of $X$ corresponding to the vertex $v\in V(G_{1})$. By symmetry, set $x(v)=x_{1}'$ for any $v\in V(K_{p_1})$, $x(v)=x_{1}$ for any $v\in V(K_{n_1-p_1-1})$, $x(v)=x_{2}$ for any $v\in V(K_{n_{2}})$, $x(v)=x_{3}$ for any $v\in V(K_{\delta -p_{1}})$, $x(v)=x_{4}$ for any $v\in V(K_{\kappa'-\delta +p_{1} })$ and $x(v)=x_{z}$ for the vertex $z$. Then 
$$X=( \underbrace{x_{1}', \ldots ,x_{1}'} _{p_{1}}, \underbrace{x_{1}, \ldots, x_{1}}_{n_{1}-p_{1}-1 }, \underbrace{x_{2}, \ldots, x_{2} }_{n_{2}}, \underbrace{x_{3}, \ldots, x_{3}}_{\delta -p_{1} }, \underbrace{x_{4}, \ldots, x_{4}}_{\kappa'-\delta +p_{1}  }, x_{z} )^{T }.$$
By $D(G_{1})X=\lambda_{1}(D(G_{1}))X $, we have 
$$\lambda_{1}(D(G_{1})) x_{z}=(\delta -p_{1} )x_{3}+2(\kappa'-\delta +p_{1})x_{4}+p_{1}x_{1}'+2(n_{1}-p_{1}-1)x_{1}+2n_{2}x_{2},$$
$$\lambda_{1}(D(G_{1})) x_{4}=(\delta -p_{1} )x_{3}+(\kappa'-\delta +p_{1}-1)x_{4}+2x_{z}+p_{1}x_{1}'+(n_{1}-p_{1}-1)x_{1}+n_{2}x_{2},$$
$$\lambda_{1}(D(G_{1})) x_{3}=(\delta -p_{1}-1 )x_{3}+(\kappa'-\delta +p_{1})x_{4}+x_{z}+p_{1}x_{1}'+(n_{1}-p_{1}-1)x_{1}+n_{2}x_{2},$$
$$\lambda_{1}(D(G_{1})) x_{2}=(\delta -p_{1} )x_{3}+(\kappa'-\delta +p_{1})x_{4}+2x_{z}+2p_{1}x_{1}'+2(n_{1}-p_{1}-1)x_{1}+(n_{2}-1)x_{2},$$
$$\lambda_{1}(D(G_{1})) x_{1}'=(\delta -p_{1} )x_{3}+(\kappa'-\delta +p_{1})x_{4}+x_{z}+(p_{1}-1)x_{1}'+(n_{1}-p_{1}-1)x_{1}+2n_{2}x_{2},$$
from which we get
$$x_{z}=\frac{(\lambda _{1}(D(G_{1}))+\kappa '-\delta +p_{1}+1)x_{4}+(n_{1}-p_{1}-1)x_{1}+n_{2}x_{2}}{\lambda _{1}(D(G_{1}))+2}>x_4$$ 
and
$$x_{3}-(x_{z}-x_{4})=\frac{(\delta-p_{1})x_{3}+2x_{z}+p_{1}x_{1}'}{\lambda_{1}(D(G_{1}))+1}>0.$$
Notice that $G_{1}$ contains $K_{\kappa'+p_{1}}$ as a proper subgraph, then $\lambda_{1}(D(G_{1}))>\lambda(D(K_{\kappa'+p_{1}}))=\kappa'+p_{1}-1$. Combining this with $\kappa'\geq \delta\geq p_1+p_2\geq 2$ and $p_{1}\geq1$, we obtain that $\lambda_{1}(D(G_{1}))>p_{1}$ and $\lambda_{1}(D(G_{1}))>\kappa'$. Similarly, $G_{1}$ contains $K_{\kappa'+n_{2}}$ as a proper subgraph, then $\lambda_{1}(D(G_{1}))>\lambda_{1}(D(K_{\kappa'+n_{2}}))=\kappa'+n_{2}-1>n_{2}$. Therefore,
\begin{equation} \label{equ1}
\begin{aligned}
&\lambda_{1}(D(G_{1}))+n_{2}+1-\frac{n_{2}(\kappa'-\delta+p_{1})}{\lambda_{1}(D(G_{1}))+1}\\
=&\frac{\lambda_{1}(D(G_{1}))^{2}+\lambda_{1}(D(G_{1}))+n_{2}\lambda_{1}(D(G_{1}))+n_{2}+\lambda_{1}(D(G_{1}))+1-n_{2}\kappa'+n_{2}\delta-n_{2}p_{1}}{\lambda_{1}(D(G_{1}))+1}\\
>&0
\end{aligned}
\end{equation}
and
\begin{equation} \label{equ2}
\begin{aligned}
&[\lambda_{1}(D(G_{1}))+n_{2}+1- \frac{n_{2}(\kappa '-\delta +p_{1})}{\lambda_{1}(D(G_{1}))+1}](x_{z}-x_{2})\\
=&(n_{2}-1)x_{z}+(\kappa'-\delta)x_{1}'+\frac{\kappa'-\delta+p_{1}}{\lambda_{1}(D(G_{1}))+1}x_{z}-\frac{n_{2}(\kappa'-\delta+p_{1})}{\lambda_{1}(D(G_{1}))+1}x_{z}\\
\geq &(n_{2}-1)x_{z}-\frac{(n_{2}-1)(\kappa'-\delta+p_{1})}{\lambda_{1}(D(G_{1}))+1}x_{z}~~(\text{since}~\kappa'\geq\delta)\\
=&\frac{(n_{2}-1)(\lambda_{1}(D(G_{1}))+1-\kappa'+\delta-p_{1})}{\lambda_{1}(D(G_{1}))+1}x_{z}\\
>&0~~(\text{since}~n_{2}\geq2~\text{and}~\lambda_{1}(D(G_{1}))>\kappa'+p_{1}-1).
\end{aligned}
\end{equation}
(\ref{equ1}) and (\ref{equ2}) mean that $x_z>x_2$. Notice that $D(G)-D(G_{1})$ is\\

\text{\qquad \qquad \quad ~~$ p_{1}$ ~~ $n_{1}-p_{1}-1$ ~~  $p_{2}$ ~~ $n_{2}-p_{2}$ ~~~~ $b$ ~~ ~~~$\delta-b-p_{1}$ ~ $\kappa'-\delta+p_{1}-1$ ~ 1 ~~~~ 1} 
\begin{center}
$\begin{matrix}
\quad \qquad \qquad ~p_{1}\\
~~~\quad n_{1}-p_{1}-1\\
~\quad \qquad \qquad p_{2}\\
~\quad \qquad n_{2}-p_{2} \\
~~\quad \qquad \qquad b\\
\quad~~~~ \delta-b-p_{1}\\
 \kappa'-\delta+p_{1}-1\\
~~~\quad \qquad \qquad 1\\
~~~\quad \qquad \qquad 1
\end{matrix}
\begin{bmatrix}
  O&\qquad  O&\qquad  O&\quad  O&~~\quad  O&\qquad~  O&\qquad\quad\quad   O&\qquad\quad  O&  O \\
  O&\qquad  O&\qquad  O&\quad  O&~~\quad  O&\qquad~  O&\qquad\quad\quad   O&\qquad\quad  O&  O \\
  O&\qquad  O&\qquad  O&\quad  O&~~\quad  O&\qquad~  O&\qquad\quad\quad   O&\qquad\quad  J& -J \\
  O&\qquad  O&\qquad  O&\quad  O&~~\quad  O&\qquad~  O&\qquad\quad\quad   O&\qquad\quad  J&  O \\
  O&\qquad  O&\qquad  O&\quad  O&~~\quad  O&\qquad~  O&\qquad\quad\quad   O&\qquad\quad  O&  O \\
  O&\qquad  O&\qquad  O&\quad  O&~~\quad  O&\qquad~  O&\qquad\quad\quad   O&\qquad\quad  O&  J \\
  O&\qquad  O&\qquad  O&\quad  O&~~\quad  O&\qquad~  O&\qquad\quad\quad   O&\qquad\quad  O&  O \\
  O&\qquad  O&\qquad  J&\quad  J&~~\quad  O&\qquad~  O&\qquad\quad\quad   O&\qquad\quad  0&  0 \\
  O&\qquad  O&\qquad -J&\quad  O&~~\quad  O&\qquad~  J&\qquad\quad\quad   O&\qquad\quad  0&  0
\end{bmatrix}.$
\end{center}
Combining the above inequalities, we obtain 
\begin{equation*}
\begin{aligned}
\lambda_{1}(D(G))-\lambda_{1}(D(G_{1}))&\ge X^{T} (D(G)-D(G_{1} ))X\\
&=2n_{2} x_{2} x_{4} +2(\delta-b-p_{1})x_{3}x_{z} -2p_{2} x_{2} x_{z}\\
&=2n_{2} x_{2} x_{4} +2p_{2}x_{3}x_{z} -2p_{2} x_{2} x_{z} ~~(\text{since}~\delta-b-p_{1}=p_{2}~\text{and}~n_{2}\geq p_{2})\\
& =2p_{2}[x_{3}x_{z}-x_{2}(x_{z}-x_{4})] \\
&>2p_{2}x_{3}(x_{z}-x_{2}) ~~(\text{since}~x_{z}>x_{4}~\text{and}~ x_{3}>x_{z}-x_{4})\\
& >0~~ (\text{since}~x_{z}>x_{2}).
\end{aligned}
\end{equation*}
It follows that $\lambda_{1}(D(G_{1}))<\lambda_{1}(D(G))$, a contradiction.

\begin{subcase}\label{subcase1.2}
$ p_{1}= n_{1}$ and $p_{2}= n_{2}$.     
\end{subcase}
Firstly, we can construct a new graph $G_{2}$ from $G$ through the following steps:
  
\text{\bf{Step 1:}} Remove $p_{2}$ edges between $z$ and $K_{p_{2}} $;
 
\text{\bf{Step 2:}} Add $p_{2}$ edges between $z$ and $S\setminus N_{S}[z]$;
 
\text{\bf{Step 3:}} Choose a vertex $w$ in $V(K_{p_{1}})$ and join $w$ to $V(K_{p_{2}})$.
\begin{figure}[H]
\centering
\includegraphics[width=0.8\linewidth]{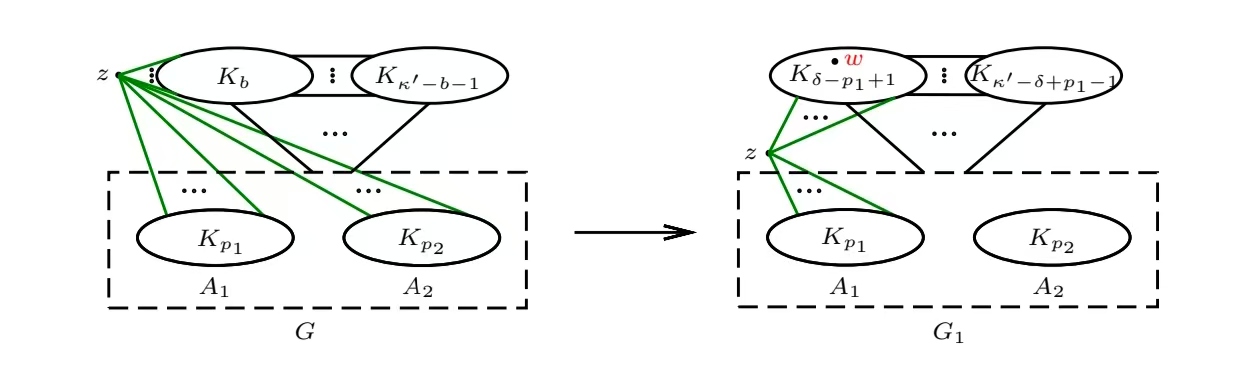}
\caption{The graphs $G$ and $G_{2}$.}
\end{figure}

Then $G_{2}-z\cong K_{\kappa '}\vee (K_{p_{1}-1}\cup K_{p_{2}})$ and $G_2\in \mathcal{G}_{n}^{\kappa', \delta}$. Let $Y$ be the Perron vector of $D(G_{2})$ and $y(v)$ denote the entry of $Y$ corresponding to the vertex $v\in V(G_{2})$. By symmetry, set $y(v)=y_{1}$ for any $v\in V(K_{p_{1}-1})$, $y(v)=y_{2}$ for any $v\in V(K_{p_{2}})$, $y(v)=y_{3}$ for any $v\in V(K_{\delta -p_{1}+1 })$, $y(v)=y_{4}$ for any $v\in V(K_{\kappa'- \delta +p_{1}-1 } )$ and $y(v)=y_{z}$ for the vertex $z$. Then 
$$Y=(\underbrace{y_{1}, \ldots, y_{1}  }_{p_{1}-1}, \underbrace{y_{2}, \ldots, y_{2}  }_{p_{2}}, \underbrace{y_{3}, \ldots, y_{3}  }_{\delta -p_{1}+1 }, \underbrace{y_{4}, \ldots, y_{4}}_{\kappa' -\delta +p_{1}-1 }, y_{z})^{T}.$$  
By $D(G_{2})Y=\lambda_{1}(D(G_{2}))Y $, we have        
$$ \lambda_{1}(D(G_{2})) y_{z}=(\delta -p_{1}+1 )y_{3}+2(\kappa'- \delta +p_{1}-1 )y_{4}+(p_{1}-1 )y_{1} +2p_{2}y_{2},$$
$$\lambda_{1}(D(G_{2})) y_{4}=(\delta -p_{1}+1 )y_{3}+(\kappa'- \delta +p_{1}-2 )y_{4}+2y_{z}+ (p_{1}-1 )y_{1}+p_{2}y_{2},$$
$$\lambda_{1}(D(G_{2})) y_{3}=(\delta -p_{1} )y_{3}+(\kappa'- \delta +p_{1}-1 )y_{4}+y_{z}+ (p_{1}-1 )y_{1} +p_{2}y_{2},$$
$$\lambda_{1}(D(G_{2})) y_{2}=(\delta -p_{1}+1 )y_{3}+(\kappa'- \delta +p_{1}-1 )y_{4}+2y_{z}+ 2(p_{1}-1 )y_{1} +(p_{2}-1)y_{2},$$ 
$$\lambda_{1}(D(G_{2})) y_{1}=(\delta -p_{1}+1 )y_{3}+(\kappa'- \delta +p_{1}-1 )y_{4}+y_{z}+ (p_{1}-2 )y_{1} +2p_{2}y_{2}.$$ 
From which we have 
$$y_{z}=y_{3}+\frac{(\kappa '-\delta +p_{1}-1)y_{4}+p_{2}y_{2}}{\lambda_{1}(D(G_{2}))+1} >y_{3}$$
and
$$y_{3}-(y_{z}-y_{3})=\frac{(\delta -p_{1}-1)y_{3}+(p_{1}-1)y_{1}}{\lambda _{1}(D(G_{2}))}>0.$$ 
Observe that $G_{2}$ contains $K_{\kappa'+p_{1}-1}$ as a proper subgraph, then $\lambda_{1}(D(G_{2}))>\lambda_{1}(D(K_{\kappa'+p_{1}-1}))=\kappa'+p_{1}-2$. Combining this with $\kappa'\geq \delta \geq p_{1}+p_{2}\geq 4$ and $p_{1} \geq 2$, we have $\lambda_{1}(D(G_{2}))>p_{1}$ and $\lambda_{1}(D(G_{2}))>\kappa'$. Similarly, $G_{2}$ also contains $K_{\kappa'+p_{2}}$ as a proper subgraph, we get $\lambda_{1}(D(G_{2}))>\lambda_{1}(K_{\kappa'+p_{2}})=\kappa'+p_{2}-1>p_{2}$. Thus,
\begin{equation} \label{equ3}
\begin{aligned}
&\lambda_{1}(D(G_{2}))+p_{2}+1-\frac{p_{2}(\kappa'-\delta+p_{1}-1)}{\lambda_{1}(D(G_{2}))+1}\\
=&\frac{\lambda_{1}(D(G_{2}))^{2}+\lambda_{1}(D(G_{2}))+p_{2}\lambda_{1}(D(G_{2}))+p_{2}+\lambda_{1}(D(G_{2}))+1-p_{2}\kappa'+p_{2}\delta-p_{1}p_{2}+p_{2}}{\lambda_{1}(D(G_{2}))+1}\\
>&0
\end{aligned}
\end{equation}
and
\begin{equation} \label{equ4}
\begin{aligned}
&[\lambda_{1}(D(G_{2}))+p_{2}+1-\frac{p_{2}(\kappa'-\delta+p_{1}-1)}{\lambda_{1}(D(G_{2}))+1}](y_{z}-y_{2})\\
=&(\kappa'-\delta)y_{1}+\frac{\kappa'-\delta+p_{1}-1}{\lambda_{1}(D(G_{2}))+1}y_{z}-\frac{p_{2}(\kappa'-\delta+p_{1}-1)}{\lambda_{1}(D(G_{2}))+1}y_{z}+(p_{2}-1)y_{z}\\
\geq& (p_{2}-1)y_{z}-\frac{(p_{2}-1)(\kappa'-\delta+p_{1}-1)}{\lambda_{1}(D(G_{2}))+1}y_{z}~~(\text{since}~\kappa'\geq\delta)\\
=& \frac{(p_{2}-1)(\lambda_{1}(D(G_{2}))+2-\kappa'+\delta-p_{1})}{\lambda_{1}(D(G_{2}))+1}\\
>&0~~(\text{since}~p_{2}\geq2~\text{and}~\lambda_{1}(D(G_{2}))>\kappa'+p_{1}-2).
\end{aligned}
\end{equation}
From (\ref{equ3}) and (\ref{equ4}), we have $y_{z}>y_{2}$. Notice that $D(G)-D(G_{2})$ is\\

\newpage
\text{\qquad \qquad \qquad \qquad    $p_{1}-1$ ~~~$p_{2}$\qquad $b$ \quad ~$\delta-p_{1}-b$ \quad ~~~1 \quad ~$\kappa'-\delta+p_{1}-1$ \quad 1}
\begin{center}
$\begin{matrix}
\qquad \quad~~~ p_{1}-1\\
\qquad \qquad \quad ~p_{2}\\
\qquad \qquad \quad ~~~b\\
\qquad \delta-p_{1}-b\\
\qquad \qquad \quad ~~~ 1\\
\kappa'-\delta+p_{1}-1\\
\qquad \qquad \quad ~~~1
\end{matrix}
\begin{bmatrix}
  O&\quad   O&\quad   O&\qquad   O&\qquad\quad   O&\qquad  \quad O&\qquad\quad  O \\
  O&\quad   O&\quad   O&\qquad   O&\qquad\quad   J&\qquad  \quad O&\qquad\quad -J \\
  O&\quad   O&\quad   O&\qquad   O&\qquad\quad   O&\qquad  \quad O&\qquad\quad  O \\
  O&\quad   O&\quad   O&\qquad   O&\qquad\quad   O&\qquad  \quad O&\qquad\quad  J \\
  O&\quad   J&\quad   O&\qquad   O&\qquad\quad   0&\qquad  \quad O&\qquad\quad  0 \\
  O&\quad   O&\quad   O&\qquad   O&\qquad\quad   O&\qquad  \quad O&\qquad\quad  O \\
  O&\quad  -J&\quad   O&\qquad   J&\qquad\quad   0&\qquad  \quad O&\qquad\quad  0
\end{bmatrix}.$
\end{center}
Then
\begin{equation*}
\begin{aligned}
\lambda_{1}(D(G))-\lambda_{1}(D(G_{2}))&\ge Y^{T} (D(G)-D(G_{2} ))Y\\
& =2p_{2}y_{2}y_{3}+2(\delta-b-p_{1})y_{3}y_{z}-2p_{2}y_{2}y_{z} \\
& =2p_{2}y_{2}y_{3}+2p_{2}y_{3}y_{z}-2p_{2}y_{2}y_{z}~~(\text{since}~\delta-b-p_{1}=p_{2}) \\
& =2p_{2}[y_{z}y_{3}-y_{2}(y_{z}-y_{3})] \\
& > 2p_{2}y_{3}(y_{z}-y_{2})~~(\text{since}~y_{z}>y_{3}~\text{and}~y_{3}>y_{z}-y_{3})\\
& >0~~ (\text{since}~y_{z}>y_{2}).
\end{aligned}
\end{equation*}
Thus, $\lambda_{1}(D(G_{2}))<\lambda_{1}(D(G))$, which contradicts the minimality of $\lambda_{1}(D(G))$. This implies that $u\notin S$.
    
By Claim \ref{claim1.1}, $z\notin S$, and so $z\in V(A_{i})$ for some $1\leq i\leq l$. Without loss of generality, we may assume that $z\in V(A_{1})$.
       
\begin{claim}\label{claim1.2}
$b=\delta-1$ .    
\end{claim}   
     
Otherwise, $b=\delta$ or $b<\delta-1$. If $b=\delta$, then $A_{1}$ is a trivial component and $V(A_{1})=\{z\}$. It means that there must be at least two non-trivial components of $A_{2}, A_{3}, \ldots, A_{l}$ and $l\geq3$. Without loss of generality, we suppose that $A_{2}$ and $A_{3}$ are two non-trivial components. By the minimality of $\lambda_{1}(D(G))$, we obtain that $G-z=K_{\kappa '} \vee (K_{n_{2}}\cup K_{n_{3}}\cup \cdots \cup K_{n_{l}})$. We first assert that $l=3$. Otherwise, we can construct a new graph $H_{1}$ obtained from $G$ by joining $K_{n_{3}}, K_{n_{4}},  \ldots , K_{n_{l}}$. Obviously, $H_1\in \mathcal{G}_{n}^{\kappa', \delta}$. According to Lemma \ref{lem1.1}, we have $\lambda_{1}(D(H_{1}))<\lambda_{1}(D(G))$, which contradicts the minimality of $\lambda_{1}(D(G))$. 
Let $G'$ be the graph obtained from $K_{1}\cup (K_{\kappa '}\vee (K_{n_{2}}\cup K_{n_{3}}))$ by add $\delta$ edges between the isolated vertex $K_{1}$ and $K_{\kappa '}$. Hence,
$$\lambda_{1} (D(G'))\leq\lambda_{1} (D(G)),$$
with equality if and only if $G\cong G'$. 
Without loss of generality, we assume that $n_{2}\leq n_{3}$. Notice that $n_{2}\geq2$ because $K_{n_{2}}$ is a non-trivial component. Furthermore, we can get that $n_2=2$. Otherwise, 
let $G''$ be the graph obtained from $K_{1}\cup (K_{\kappa '}\vee (K_{2}\cup K_{n-\kappa '-3}))$ by adding $\delta$ edges between $K_{1}$ and $K_{\kappa'}$. Obviously, $G'\in \mathcal{G}_{n}^{\kappa', \delta}$ and $G''\in \mathcal{G}_{n}^{\kappa', \delta}$.
Let $X$ be the Perron vector of $D(G'')$ and $x(v)$ denote the entry of $X$ corresponding to the vertex $v\in V(G'')$. By symmetry, we say $x(v)=x_{1}$ for the vertex $v\in V(K_2)$, $x(v)=x_{2}$ for the vertex $v\in V(K_{n-\kappa'-3})$, $x(v)=x_{3}$ for the vertex $v\in V(K_\delta)$, $x(v)=x_{4}$ for the vertex $v\in V(K_{\kappa'-\delta})$ and $x(v)=x_{z}$ for the vertex $z$. Note that $n-\kappa'-3=n_{2}+n_3-2$. Then
$$X=(x_1,x_1,\underbrace{x_2,\ldots,x_2 }_{n_2+n_3-2},\underbrace{x_3,\ldots,x_3 }_{\delta },\underbrace{x_4,\ldots,x_4 }_{\kappa '-\delta },x_{z}).$$
By $D(G'')X=\lambda_{1}(D(G''))X $, we have 
$$\lambda_{1}(D(G''))x_{2}=\delta x_{3}+(\kappa'-\delta)x_{4}+2x_{z}+4x_{1}+(n_{2}+n_{3}-3)x_{2},$$
$$\lambda_{1}(D(G''))x_{1}=\delta x_{3}+(\kappa'-\delta)x_{4}+2x_{z}+x_{1}+2(n_{2}+n_{3}-2)x_{2},$$
$$\lambda_{1}(D(G''))x_{z}=\delta x_{3}+2(\kappa'-\delta)x_{4}+4x_{1}+2(n_{2}+n_{3}-2)x_{2}.$$
Since $\kappa' \geq \delta$ and $n_{3}\geq n_{2}\geq3$, $x_{z}-x_{2}=\frac{(\kappa'-\delta)x_{4}+(n_{2}+n_{3}-3)x_{2}}{\lambda_{1}(D(G''))+2}>0$ and $x_{1}-x_{2}=\frac{(n_{2}+n_{3}-4)x_{2}}{\lambda_{1}(D(G''))+3}>0$. Then $x_{z}>x_{2}$ and $x_{1}>x_{2}$. Notice that $D(G')-D(G'')$ is\\

\text{\qquad \qquad \qquad\qquad\qquad\quad~~~~2~\quad~~$n_{2}-2$~\quad$n_{3}$~\quad~~$\delta$\quad$\kappa'-\delta$~~~~~1}
\begin{center}
$\begin{matrix}
\quad~~~~2\\
n_{2}-2\\
\quad~~~n_{3}\\
\quad~~~~\delta\\
\kappa'-\delta\\
\quad~~~~1
\end{matrix}$
$\begin{bmatrix}
  O&\quad -J&\quad  O&\quad  O&\quad  O&\quad  O \\
 -J&\quad  O&\quad  J&\quad  O&\quad  O&\quad  O \\
  O&\quad  J&\quad  O&\quad  O&\quad  O&\quad  O \\
  O&\quad  O&\quad  O&\quad  O&\quad  O&\quad  O \\
  O&\quad  O&\quad  O&\quad  O&\quad  O&\quad  O \\
  O&\quad  O&\quad  O&\quad  O&\quad  O&\quad  0
\end{bmatrix}.$
\end{center} 
Combining the above inequalities, we obtain 
\begin{equation*}
\begin{aligned}
&\lambda_{1}(D(G'))-\lambda_{1}(D(G''))\\
\ge& X^{T} (D(G')-D(G''))X\\
=&2n_{3}(n_{2}-2)x_{2}^2-4(n_{2}-2)x_{1}x_{2}\\
=&2(n_{2}-2)x_{2}(n_{3}x_{2}-2x_{1})\\
>&(n_{3}-2)\delta x_{3}+(n_{3}-2)(\kappa'-\delta)x_{4}+(2n_{3}-4)x_{z}+(4n_{3}-2)x_{1}+(n_{2}n_{3}+n_{3}^{2}-7n_{3}-4n_{2}+8)x_{2}\\
&(\text{since}~n_{3}\geq3)\\
>&(6n_{3}-6)x_{2}+(n_{2}n_{3}+n_{3}^{2}-7n_{3}-4n_{2}+8)x_{2}~~(\text{since}~\kappa'\geq\delta, x_{z}>x_{2}~\text{and}~x_{1}>x_{2})\\
=&[n_{3}(n_{3}-1)+n_{2}(n_{3}-4)+2]x_{2}\\
\geq & [n_{2}(n_{3}-1)+n_{2}(n_{3}-4)+2]x_{2}~~(\text{since}~n_{3}\geq n_{2})\\
=&[n_{2}(2n_{3}-5)+2]x_{2}\\
>&0~~(\text{since}~n_{3}\geq3).
\end{aligned}
\end{equation*}
Therefore, $\lambda_{1}(D(G''))<\lambda_{1}(D(G'))$, a contradiction. By the minimality of $\lambda_{1}(D(G))$, we infer that $G-z\cong K_{\kappa '}\vee (K_{2}\cup K_{n-\kappa '-3})$, i.e., $(n_{2}, n_{3})=(2, n-\kappa'-3)$. Denote $V(K_{2})=\{w_{1}, w_{2}\}$ and $N_{S}(z)= \{ u_{1}, u_{2}, \ldots ,u_{\delta } \}$. Next, we can construct a new graph $G_{3}$ from $G$ by the following steps:

\text{\bf{Step 1:}} Remove $zu_{\delta}$ and $w_{1}w_{2}$;

\text{\bf{Step 2:}} Add $zw_{1}$ and edges between $w_{2}$ and $V(K_{n_{3}})$.
 \begin{figure}[H]
\centering
\includegraphics[width=0.8\linewidth]{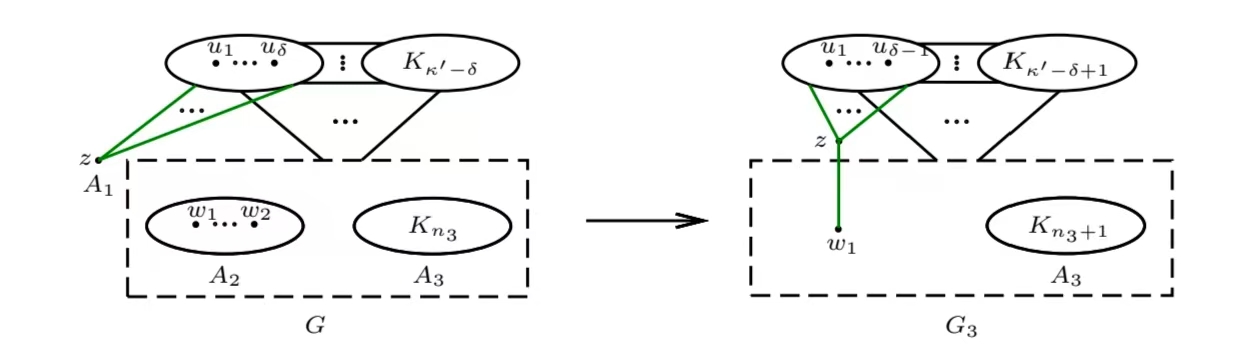}
\caption{The graphs $G$ and $G_{3}$.}
\end{figure}

Then $G_{3}-z=K_{\kappa'}\vee (K_{1}\cup K_{n_{3}+1})$ and $G_3\in \mathcal{G}_{n}^{\kappa', \delta}$, where $n_{3}=n-\kappa'-3$. Let $X$ be the Perron vector of $D(G_{3})$ and $x(v)$ denote the entry of $X$ corresponding to the vertex $v\in V(G_{3})$. By symmetry, set $x(v)=x_{w_{1}}$ for the vertex $w_{1}$,  $x(v)=x_{2}$ for any $v\in V(K_{n_{3}+1})$, $x(v)=x_{3}$ for any $v\in V(K_{\delta-1})$, $x(v)=x_{4}$ for any $v\in V(K_{\kappa'-\delta+1})$ and $x(v)=x_{z}$ for the vertex $z$, then
$$X=(x_{w_{1}}, \underbrace{x_{2},\ldots, x_{2}  }_{n_{3}+1 }, \underbrace{x_{3}, \ldots, x_{3}}_{\delta -1}, \underbrace{x_{4}, \ldots, x_{4}}_{\kappa ' -\delta +1 }, x_{z})^{T}.$$ 
By $D(G_{3})X=\lambda_{1}(D(G_{3}))X$, we have\\
$$\lambda_{1}(D(G_{3})) x_{w_{1}}=(\delta -1)x_{3}+(\kappa'-\delta +1 )x_{4}+x_{z}+2(n_{3}+1 )x_{2},$$ 
$$\lambda_{1}(D(G_{3})) x_{2}=(\delta -1)x_{3}+(\kappa'-\delta +1 )x_{4}+2x_{z}+2x_{w_{1}}+ n_{3}x_{2},$$ 
$$\lambda_{1}(D(G_{3})) x_{z}=(\delta -1)x_{3}+2(\kappa'-\delta +1 )x_{4}+x_{w_{1}}+ 2(n_{3}+1)x_{2},$$ 
$$\lambda_{1}(D(G_{3})) x_{4}=(\delta -1)x_{3}+(\kappa'-\delta )x_{4}+2x_{z}+ x_{w_{1}}+ (n_{3}+1)x_{2}.$$ 
From which we have
\begin{equation} \label{equ5}
\begin{aligned}
2x_{2}-x_{w_{1}}=\frac{(\delta -1)x_{3}+(\kappa '-\delta +1)x_{4}+3x_{z}+3x_{w_{1}}}{(\lambda_{1}(D(G_{3})) +1)}
> 0
\end{aligned}
\end{equation}
and 
\begin{equation} \label{equ6}
\begin{aligned}
&x_{w_{1}}-x_{4}\\
=&\frac{(n_3+1)x_2-x_z}{\lambda_1(D(G_3))+1}\\
=&\frac{\lambda_1(D(G_3))[(n_3+1)x_2-x_z]}{\lambda_1(D(G_3))(\lambda_1(D(G_3))+1)}\\
=&\frac{n_3(\delta-1)x_3+(n_3-1)(\kappa'-\delta+1)x_4+2(n_3+1)x_z+(2n_3+1)x_{w_1}+(n_3+1)(n_3-2)x_2}
{\lambda_1(D(G_3))(\lambda_1(D(G_3))+1)}\\
>&0~~(\text{since}~\delta>1,~n_3\geq2~\text{and}~\kappa'\geq \delta).
\end{aligned}
\end{equation}
Notice that $D(G)-D(G_{3})$ is\\

\text{\qquad \qquad \qquad \qquad \qquad \quad~~ $1 \quad ~~1~~~~~n_{3}~~~~\delta-1~~~~~~1~~~~~\kappa'-\delta~~~1$}  
\begin{center}
$\begin{matrix}
\quad ~~~~1\\
\quad ~~~~ 1\\
\quad ~~~ n_{3}\\
~\delta-1\\
\quad ~~~~1\\
\kappa'-\delta\\
\quad ~~~~1
\end{matrix}
\begin{bmatrix}
  0& -1&  O&&  O&&  0&&  O&~  1 \\
 -1&  0&  J&&  O&&  0&&  O&~  0 \\
  O&  J&  O&&  O&&  O&&  O&~  O \\
  O&  O&  O&&  O&&  O&&  O&~  O \\
  0&  0&  O&&  O&&  0&&  O&~ -1 \\
  O&  O&  O&&  O&&  O&&  O&~  O \\
  1&  0&  O&&  O&& -1&&  O&~  0 
\end{bmatrix}.$
\end{center}       
Therefore,
\begin{equation*}     
 \begin{aligned}
\lambda_{1} (D(G))-\lambda_{1} (D(G_{3}) )&\ge X^{T} (D(G)-D(G_{3} ))X\\
&=2n_{3}x_{2}^{2}+2x_{w_{1}}x_{z}-2x_{2}x_{w_{1}}-2x_{4}x_{z}\\
&=2x_{2}(n_{3}x_{2}-x_{w_{1}})+2x_{z}(x_{w_{1}}-x_{4})\\
&\ge 2x_{2}(2x_{2}-x_{w_{1}})+2x_{z}(x_{w_{1}}-x_{4})~~(\text{since}~~n_{3}\geq2)\\
&>0~~(\text{by}~(\ref{equ5})~\text{and}~(\ref{equ6})). 
\end{aligned}    
\end{equation*} 
Thus, $\lambda_{1}(D(G_{3}))<\lambda_{1}(D(G)$), which contradicts the minimality of $\lambda_{1}(D(G))$. This demonstrates that $b\neq\delta$.     
      
If $b<\delta-1$, then $d_{A_{1}}(z)=\delta -b\geq2$ and $A_{1}$ is a non-trivial component. Without loss of generality, we assume that $A_{2}$ is also a non-trivial component. By the minimality of $\lambda_{1}(D(G))$, we can obtain that $G-z\cong K_{\kappa '}\vee (K_{n_{1}-1}\cup K_{n_{2}}\cup \cdots \cup K_{n_{l}})$. We first assert that $l=2$. Otherwise, we construct a new graph $H_{2}$ obtained from $G$ by joining $K_{n_{2}}, K_{n_{3}},  \ldots , K_{n_{l}}$. Obviously, $H_2\in \mathcal{G}_{n}^{\kappa', \delta}$. According to Lemma \ref{lem1.1}, we obtain that $\lambda_{1}(D(H_{2}))<\lambda_{1}(D(G))$, which is a contradiction. Let $N_{S}(z)=\{u_{1}, u_{2}, \ldots ,u_{b}\}$, $S\backslash N_{S}(z)=\{u_{b+1}, u_{b+2}, \ldots, u_{\kappa'}\}$ and $N_{A_{1}}(z)= \{ w_{1}, w_{2}, \ldots ,w_{\delta -b}\}$. Next, we can construct a new graph $G_{4}$ from $G$ by the following steps:
  
\text{\bf{Step 1:}} Remove $\delta-b-1$ edges between $z$ and $N_{A_1}(z)$ and retain $zw_1$;
  
\text{\bf{Step 2:}} Add $\delta-b-1$ edges between $z$ and $S\setminus N_{S}(z)$;
 
\text{\bf{Step 3:}} Remove $n_1-2$ edges between $w_1$ and $V(A_1)\backslash \{z,w_1\}$;
 
\text{\bf{Step 4:}} Add all edges between $V(A_1)\backslash \{z,w_1\}$ and $V(A_2)$.
\begin{figure}[H]
\centering
\includegraphics[width=\linewidth]{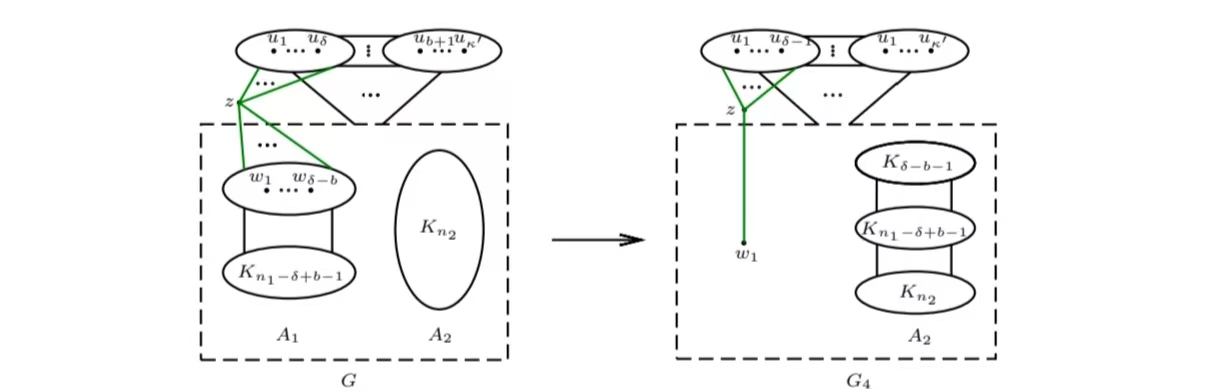}
\caption{The graphs $G$ and $G_{4}$.}
\end{figure}

Obviously, $G_4\in \mathcal{G}_{n}^{\kappa', \delta}$.
Let $Y$ be the Perron vector of $D(G_{4})$ and $y(v)$ denote the entry of $Y$ corresponding to the vertex $v\in V(G_{4})$. By symmetry, set $y(v)=y_{w_{1}}$ for the vertex $w_{1}$,  $y(v)=y_{2}$ for any $v\in V(K_{n_{1}+n_{2}-2})$, $y(v)=y_{3}$ for any $v\in V(K_{\delta-1})$, $y(v)=y_{4}$ for any $v\in V(K_{\kappa'-\delta+1})$ and $y(v)=y_{z}$ for the vertex $z$, then 
$$Y=(y_{w_{1}}, \underbrace{y_{2}, \ldots, y_{2} }_{n_{1}+n_{2}-2}, \underbrace{y_{3}, \ldots, y_{3}}_{\delta -1}, \underbrace{y_{4},\ldots,y_{4}}_{\kappa'-\delta +1}, y_{z})^{T}.$$ 
By $D(G_{4})Y=\lambda_{1}(D(G_{4}))Y $, we have\\
$$\lambda_{1}(D(G_{4})) y_{2}=(\delta -1)y_{3}+(\kappa'-\delta +1 )y_{4} +2y_{z} +2y_{w_{1}}+(n_{1}+n_{2}-3)y_{2},$$
$$\lambda_{1}(D(G_{4})) y_{w_{1}}=(\delta -1)y_{3}+(\kappa'-\delta +1 )y_{4} +y_{z} +2(n_{1}+n_{2}-2)y_{2},$$
$$\lambda_{1}(D(G_{4})) y_{3}=(\delta -2)y_{3}+(\kappa'-\delta +1 )y_{4} +y_{z} +y_{w_{1}}+(n_{1}+n_{2}-2)y_{2},$$
$$\lambda_{1}(D(G_{4})) y_{z}=(\delta -1)y_{3}+2(\kappa'-\delta +1 )y_{4} +y_{w_{1}} +2(n_{1}+n_{2}-2   )y_{2}.$$
From which we get

\begin{equation} \label{equ7}
\begin{aligned}
2y_2-y_z=\frac{(\delta -1)y_{3}+3y_{z}+3y_{w_{1}} }{\lambda_{1}(D(G_{4}))+1}>0,
\end{aligned}
\end{equation}

\begin{equation} \label{equ8}
\begin{aligned}
y_{w_1}-y_3=\frac{(n_{1}+n_{2}-2)y_{2}}{\lambda_{1}(D(G_{4})) +1}>0,
\end{aligned}
\end{equation}

\begin{equation} \label{equ9}
\begin{aligned}
y_z-y_2=\frac{(\kappa '-\delta +1)(\lambda_{1}(D(G_{4}))+n_{1}+n_{2}-2)y_4}{(\lambda_{1}(D(G_{4})) +1)(\lambda_{1}(D(G_{4}))+n_{1}+n_{2}-1)}
+\frac{(n_{1}+n_{2}-4)y_{w_{1}}}{(\lambda_{1}(D(G_{4}))+n_{1}+n_{2}-1)}>0.
\end{aligned}
\end{equation}
Notice that $D(G)-D(G_{4})$ is\\
\newpage
\text{\qquad \qquad \qquad  ~ $1$~~~~~~$\delta-b-1$~~~~$n_{1}-\delta+b-1$~~~~~$~n_{2}$~~~~~$b$~~~~~~$\delta-b-1$~~$\kappa'-\delta+1$~~~~~~1} 
\begin{center}
$\begin{matrix}
\qquad \qquad ~~~~~1\\
\quad ~~~~\delta-b-1\\
n_{1}-\delta+b-1\\
\qquad \qquad ~~~n_{2}\\
\qquad \qquad ~~~~b\\
\quad ~~~~\delta-b-1\\
\quad ~~\kappa'-\delta+1\\
\qquad \qquad ~~~~~1
\end{matrix}$
$\begin{bmatrix}
  0&\quad~~ -J&\qquad\quad~ -J&\qquad\quad   O&~~  O&\quad\quad   O&\qquad~~~  O&\qquad  0 \\
 -J&\quad~~  O&\qquad\quad~  O&\qquad\quad   J&~~  O&\quad\quad   O&\qquad~~~  O&\qquad -J \\
 -J&\quad~~  O&\qquad\quad~  O&\qquad\quad   J&~~  O&\quad\quad   O&\qquad~~~  O&\qquad  O \\
  O&\quad~~  J&\qquad\quad~  J&\qquad\quad   O&~~  O&\quad\quad   O&\qquad~~~  O&\qquad  O \\
  O&\quad~~  O&\qquad\quad~  O&\qquad\quad   O&~~  O&\quad\quad   O&\qquad~~~  O&\qquad  O \\
  O&\quad~~  O&\qquad\quad~  O&\qquad\quad   O&~~  O&\quad\quad   O&\qquad~~~  O&\qquad  J \\
  O&\quad~~  O&\qquad\quad~  O&\qquad\quad   O&~~  O&\quad\quad   O&\qquad~~~  O&\qquad  O \\
  0&\quad~~ -J&\qquad\quad~  O&\qquad\quad   O&~~  O&\quad\quad   J&\qquad~~~  O&\qquad  0
\end{bmatrix}.$ 
\end{center}       
Then,      
\begin{equation*}
\begin{aligned}
&\lambda_{1} (D(G))-\lambda_{1} (D(G_{4}) )\\
\ge& Y^{T} (D(G)-D(G_{4}))Y\\
=& 2n_{2}(n_{1}-2)y_{2}^{2}+2(\delta -1-b)y_{z}y_{3}-2(\delta -1-b)y_{2}y_{z}-2(n_{1}-2)y_{2}y_{w_{1}}\\ 
=& 2(n_{1}-2)y_{2}(n_{2}y_{2}-y_{w_{1}} )+2(\delta -1-b)y_{z}(y_{3}-y_{2}y_{z}) \\
\ge& 2(n_{1}-2)y_{2}(2y_{2}-y_{w_{1}} )+2(\delta -1-b)y_{z}(y_{3}-y_{2}y_{z})~~(\text{since}~~n_{2}\ge 2)\\
\ge& 2(\delta -1-b)y_{2}(2y_{2}-y_{w_{1}})+2(\delta -1-b)y_{z}(y_{3}-y_{2}y_{z})
~~(\text{since}~n_{1}\!-\!2\ge\delta \!-\!1\!-b)\\
=& 2(\delta -1-b)y_{2}\left[\frac{(\delta -1)y_{3} }{\lambda_{1}(D(G_{4})) }+\frac{(\kappa'-\delta +1 )y_{4} }{\lambda_{1}(D(G_{4})) } +\frac{3y_{z} }{\lambda_{1}(D(G_{4})) }+\frac{4y_{w_{1}} }{\lambda_{1}(D(G_{4})) }-\frac{2y_{2} }{\lambda_{1}(D(G_{4})) }\right]\\ 
&+2(\delta -1-b)y_{z}\left[ \frac{y_{2}  }{\lambda_{1}(D(G_{4})) }-\frac{y_{3} }{\lambda_{1}(D(G_{4}))} -\frac{y_{z} }{\lambda_{1}(D(G_{4})) }-\frac{y_{w_{1}} }{\lambda_{1}(D(G_{4})) }\right]\\
\end{aligned}
\end{equation*}
\begin{equation*}
\begin{aligned}
=& 2(\delta -1-b)\frac{(\delta -1)y_{2}y_{3}+(\kappa'-\delta +1 )y_{2}y_{4}+4y_{2}y_{z}+4y_{2}y_{w_{1}}-2y_{2}^{2}-y_{z}y_{3}-y_{z}^{2}-y_{z}y_{w_{1}}}{\lambda_{1}(D(G_{4}))}\\ 
>& 2(\delta -1-b)\frac{y_z(2y_2-y_z)+2y_2(y_z-y_2)+y_{w_1}(2y_2-y_z)+(2y_2y_{w_1}-y_zy_3)}{\lambda_{1}(D(G_{4})) }\\
&(\text{since}~\delta>1~\text{and}~\kappa'\geq\delta)\\
>& 0~~(\text{by}(\ref{equ7})-(\ref{equ9})).        
\end{aligned}
\end{equation*}
Thus, we get $\lambda_{1}(D(G_{4}))<\lambda_{1}(D(G))$, which gives a contradiction. Therefore, we conclude that $b=\delta-1$, completing the proof of Claim \ref {claim1.2}.

Without loss of generality, we still suppose that $z\in V(A_{1})$. According to Claim \ref {claim1.2}, we get $d_{A_{1}}(z)=1$. For convenience, we denote $N_{A_{1}}(z)=\{w\}$. By the minimality of $\lambda_{1}(D(G))$, we can infer that $G-z\cong K_{\kappa '}\vee (K_{n_{1}-1}\cup K_{n_{2}})$. That is, $G-S$ has exactly two non-trivial components $A_{1}$ and $A_{2}$. Recall that $V(A_{i})=n_{i}$, where $n_{i}\geq2$ for $i=1,2$. 

\begin{claim}\label{claim1.3}
$n_{1}=2$ .    
\end{claim}

Otherwise, $n_{1}\geq3$. Assume that $E_{1}''=\{v_{i}u_{i}\mid v_{i}\in V(A_{1})\backslash\{z,w\}, u_{i}\in V(A_{2})\}$ and $E_{2}''=\{wv_{i}\mid v_{i}\in V(A_{1})\backslash\{z,w\}\}$. Let $G_5=G+E_1''-E_2''$. Obviously, $G_{5}-z\cong K_{\kappa'}\vee(K_{1}\cup K_{n-\kappa'-2})$ and $G_5\in \mathcal{G}_{n}^{\kappa', \delta}$. Note that $n-\kappa'=n_{1}+n_{2}$. 
Let $X$ be the Perron vector of $D(G_{5})$ and let $x(v)$ denote the entry of $X$ corresponding to the vertex $v\in V(G_{5})$. By symmetry, set $x(v)=x_{w}$ for the vertex $w$, $x(v)=x_{2}$ for any $v\in V(K_{n_{1}+n_{2}-2})$, $x(v)=x_{3}$ for any $v\in V(K_{\delta-1})$, $x(v)=x_{4}$ for any $v\in V(K_{\kappa'-\delta+1})$ and $x(v)=x_{z}$ for the vertex $z$, then   
$$X=(x_{w},\underbrace{x_{2},\ldots ,x_{2}}_{n_{1}+n_{2}-2},\underbrace{x_{3},\ldots x_{3}}_{\delta -1},\underbrace{x_{4},\ldots ,x_{4}}_{\kappa '-\delta +1},x_{z})^{T}.$$
By $D(G_{5})X=\lambda_{1}(D(G_{5}))X$, we have
 
$$\lambda _{1}(D(G_{5}))x_{2}=(\delta -1)x_{3}+(\kappa '-\delta +1)x_{4}+2x_{z}+2x_{w}+(n_{1}+n_{2}-3)x_{2},$$
$$\lambda _{1}(D(G_{5}))x_{w}=(\delta -1)x_{3}+(\kappa '-\delta +1)x_{4}+x_{z}+2(n_{1}+n_{2}-2)x_{2},$$
from which we get
\begin{equation*}
\begin{aligned}
2x_{2}-x_{w}=\frac{(\delta-1)x_{3}+(\kappa'-\delta+1)x_{4}+3x_{z}+3x_{w}}{\lambda_{1}(D(G_{5}))+1}>0.
\end{aligned}
\end{equation*}
Notice that $D(G)-D(G_{5})$ is\\

\text{\qquad \qquad \qquad \qquad \qquad \qquad ~~1 ~~~$n_{1}-2$ ~~ $n_{2}$~~~~$\delta-1$~~~$\kappa'-\delta+1$ ~ 1}
\begin{center}
$\begin{matrix}
\qquad \quad~ 1\\ 
\quad~ n_{1}-2\\
\qquad \quad n_{2}\\
\quad ~~~ \delta-1\\
\kappa'-\delta+1\\
\qquad \quad~~ 1
\end{matrix}$ 
$\begin{bmatrix}
  0& -J&\quad  O&&  O&&\quad  O&&\quad  0 \\
 -J&  O&\quad  J&&  O&&\quad  O&&\quad  O \\
  O&  J&\quad  O&&  O&&\quad  O&&\quad  O \\
  O&  O&\quad  O&&  O&&\quad  O&&\quad  O \\
  O&  O&\quad  O&&  O&&\quad  O&&\quad  O \\
  0&  O&\quad  O&&  O&&\quad  O&&\quad  0
\end{bmatrix}.$
\end{center}
Thus, 
\begin{align*}
\lambda_{1} (D(G))-\lambda_{1} (D(G_{5}))&\ge X^{T} (D(G)-D(G_{5} ))X\\
&=2n_{2}(n_{1}-2)x_{2}^{2}-2(n_{1}-2)x_{2}x_{w}\\
&=2x_{2}(n_{1}-2)(n_{2}x_{2}-x_{w})\\
&\geq2(n_{1}-2)x_{2}(2x_{2}-x_{w})~~(\text{since}~n_{2}\geq2)\\
&>0~~(\text{since}~~n_{1}\geq 3 ~\text{and}~ 2x_{2}> x_{w}).
\end{align*} 
It follows that $\lambda_{1}(D(G_{5}))<\lambda_{1}(D(G))$, a contradiction. Hence $n_{1}=2$. 
Then we can deduce that $G\cong G_{n}^{\kappa',\delta}$, as desired.

\begin{case}\label{case1.2}   
$\kappa '\leq \delta -1$.
\end{case}
By the minimality of $\lambda_1(D(G))$, we can assert that $l=2$. Otherwise, $l\geq3$. Without loss of generality, we suppose that $A_{1}$ is one of the non-trivial components. Let $G'=K_{\kappa '}\vee (K_{n_{1}}\cup K_{n_{2}'})$, where $n_{2}'=\sum_{i=2}^{l}n_{i}$. According to Lemma \ref{lem1.1}, it is easy to see that $\lambda_{1}(D(G'))<\lambda_{1}(D(G))$ , a contradiction. This means that $l=2$. Therefore, $G-S$ has exactly two non-trivial components $A_{1}$ and $A_{2}$. Furthermore,
$$\lambda _{1}(D(K_{\kappa '}\vee (K_{n_{1}}\cup K_{n_{2}})))\leq\lambda _{1}(D(G)),$$
with equality if and only if $G\cong K_{\kappa '}\vee (K_{n_{1}}\cup K_{n_{2}})$.
Without loss of generality, we may assume that $n_{1}\geq n_{2}$. Obviously, $n_{2}\geq \delta-\kappa'+1$ because the minimum degree of $G$ is $\delta$. Combining this with Lemma \ref{lem1.2}, we get
$$\lambda _{1}(D(K_{\kappa '}\vee (K_{n-\delta-1}\cup K_{\delta-\kappa'+1})))\leq\lambda _{1}(D(K_{\kappa '}\vee (K_{n_{1}}\cup K_{n_{2}}))),$$
with equality if and only if $(n_{1}, n_{2})=(n-\delta-1, \delta-\kappa'+1)$. By the minimality of $\lambda _{1}(D(G))$, we infer that $G\cong K_{\kappa '}\vee (K_{n-\delta-1}\cup K_{\delta-\kappa'+1})$.

This completes the proof.
\begin{flushright}
$\square$
\end{flushright}

\section{Proof of Theorem \ref{thm1.3}}\label{3}
\begin{lemma}\cite{H.Q. Lin}\label{lem1.3}
Let $\overrightarrow{G}$ be a strongly connected digraph with $u,v\in V(\overrightarrow{G})$ and $uv \notin E(\overrightarrow{G})$. Then $\lambda_{1}(D(\overrightarrow{G}))>\lambda_{1}(D(\overrightarrow{G}+uv))$.
\end{lemma}

\begin{lemma}\cite{J.A. Bondy}\label{lem1.4}
Let $\overrightarrow{G}$ be an arbitrary strongly connected digraph with vertex connectivity $k$. Suppose that $S$ is a $k$-vertex cut of $\overrightarrow{G}$ and $\overrightarrow{G_{1}}, \ldots, \overrightarrow{G_{s}}$ are the strongly connected components of $\overrightarrow{G}-S$. Then there exists an ordering of $\overrightarrow{G}_{1}, \ldots, \overrightarrow{G}_{s}$ such that, for $1\leq i\leq s$ and $v\in V(\overrightarrow{G}_{i})$, every tail of $v$ in $\overrightarrow{G}_{1}, \ldots, \overrightarrow{G}_{i-1}$.
\end{lemma}

\begin{lemma}\label{lem1.5}
Let $f(x)=-4x^{2}+4(n-k)x+4n-4$, where $2\leq x\leq n-k-2$ and $k\geq1$. Then $f_{min}(x)=f(2)=f(n-k-2)=12n-8k-20$.
\end{lemma}

\noindent{\bf{Proof of Theorem \ref{thm1.3}}.}
Assume that $\overrightarrow{G}$ is a digraph with minimum distance spectral radius in $\overrightarrow{\mathcal{D}}_{n,k}$. Then there exists some subset $S\subseteq V(\overrightarrow{G})$ with $|S|=k$ such that $\overrightarrow{G}-S$ contains at least two non-trivial strongly connected components. Let $\overrightarrow{G}_{1},\ldots, \overrightarrow{G}_{l}$ be the $l$ strongly connected components of $\overrightarrow{G}-S$ and $\overrightarrow{G}_{i}$ be the first non-trivial strongly connected component. 
It is easy to see that $1\leq i\leq l-1$ and there must be a non-trivial strongly connected component $\overrightarrow{G}_{j}$ for $i+1\leq j\leq l$. According to Lemma \ref{lem1.3}, let $\overrightarrow{H}_{1}= \overrightarrow{G}_{1}\nabla\cdots \nabla\overrightarrow{G}_{i}$, $\overrightarrow{H}_{2}= \overrightarrow{G}_{i+1}\nabla \cdots \nabla \overrightarrow{G}_{l}$ and $|V(\overrightarrow{H}_{1})|=n_1$. We can construct a new graph $\overrightarrow{G}'$ obtained from $\overrightarrow{G}$ by adding arcs until induced subdigraph $V(\overrightarrow{H}_{1})\cup S$ and induced subdigraph $V(\overrightarrow{H}_{2})\cup S$ are complete digraphs. There exists an ordering of $\overrightarrow{H}_{1}$ and $\overrightarrow{H}_{2}$ such that every tail of $v\in V(\overrightarrow{H_{2}})$ in $\overrightarrow{H}_{1}$ by Lemma \ref{lem1.4}.
Obviously, $\overrightarrow{G}'= \overrightarrow{G}_{n}^{k,n_{1}}\in \overrightarrow{\mathcal{G}}(n,k)\subseteq \overrightarrow{\mathcal{D}}_{n,k}$. Then $\lambda_{1}(D(\overrightarrow{G}'))\leq \lambda_{1}(D(\overrightarrow{G}))$, with equality if and only if $\overrightarrow{G}'\cong \overrightarrow{G}$. Thus, the digraphs which achieve the minimum distance spectral radius among all digraphs in $\overrightarrow{\mathcal{D}}_{n,k}$ must be in $\overrightarrow{\mathcal{G}}(n,k)$. By the minimality of $\lambda_{1}(D(\overrightarrow{G}))$, we get $\overrightarrow{G}\cong \overrightarrow{G}_{n}^{k,n_{1}}$ for any integer $n_{1}\geq 2$ and $n-k-n_{1}\geq 2$. For convenience, let $n-k-n_{1}=n_{2}$. Then
\begin{center}
$D(\overrightarrow{G})$=
$\begin{bmatrix}
  J_{n_{1}\times n_{1}}-I_{n_{1}\times n_{1}}&  J_{n_{1}\times n_{2}}&                        J_{n_{1}\times k}& \\
  2J_{n_{2}\times n_{1}}&                       J_{n_{2}\times n_{2}}-I_{n_{2}\times n_{2}}&  J_{n_{2}\times k}& \\
  J_{k\times n_{1}}&                            J_{k\times n_{2}}&                            J_{k\times k}-I_{k\times k}&
\end{bmatrix}$
\end{center}
and
\begin{equation*}
\begin{aligned}
&P(\lambda(D(\overrightarrow{G})))\\
=&|\lambda(D(\overrightarrow{G}))I-D(\overrightarrow{G})|\\
=&\begin{vmatrix}
  (\lambda(D(\overrightarrow{G}))+1)I_{n_{1}\times n_{1}}-J_{n_{1}\times n_{1}}&  -J_{n_{1}\times n_{2}}&                                                             -J_{n_{1}\times k}& \\
  -2J_{n_{2}\times n_{1}}&                                                            (\lambda(D(\overrightarrow{G}))+1)I_{n_{2}\times n_{2}}-J_{n_{2}\times n_{2}}&  -J_{n_{2}\times k}& \\
  -J_{k\times n_{1}}&                                                                 -J_{k\times n_{2}}&                                                                 (\lambda(D(\overrightarrow{G}))+1)I_{k\times k}-J_{k\times k}&
\end{vmatrix}\\
=&(\lambda(D(\overrightarrow{G}))+1)^{n-3}
\begin{vmatrix}
  \lambda(D(\overrightarrow{G}))-n_{1}+1&  -n_{2}&                                      -k& \\
  -2n_{1}&                                     \lambda(D(\overrightarrow{G}))-n_{2}+1&  -k& \\
  -n_{1}&                                      -n_{2}&                                      \lambda(D(\overrightarrow{G}))-k+1&
\end{vmatrix}\\
=&(\lambda(D(\overrightarrow{G}))+1)^{n-2}
\begin{vmatrix}
  \lambda(D(\overrightarrow{G}))-n_{1}-k+1&  -n_{2}&                                      \\
  -2n_{1}-k&                                     \lambda(D(\overrightarrow{G}))-n_{2}+1&   \\                                                                     
\end{vmatrix}\\
=&(\lambda(D(\overrightarrow{G}))+1)^{n-2}[\lambda(D(\overrightarrow{G}))^{2}-(n-2) \lambda(D(\overrightarrow{G}))-n_{1}n_{2}-n+1]\\
=&(\lambda(D(\overrightarrow{G}))+1)^{n-2}[\lambda(D(\overrightarrow{G}))^{2}-(n-2)\lambda(D(\overrightarrow{G}))-n_{1}(n-k-n_{1})-n+1].
\end{aligned}
\end{equation*}
It shows that $\lambda_{1}(D(\overrightarrow{G}))$ is equal to the largest root of the equation
$$x^{2}-(n-2)x-n_{1}(n-k-n_{1})-n+1=0,$$
where $2\leq n_{1}\leq n-k-2$. Then according to Lemma \ref{lem1.5}, we have 
$$\lambda_{1}(D(\overrightarrow{G}))=\frac{n-2+\sqrt{(n-2)^2-4[n_{1}^{2}-(n-k)n_{1}-n+1]}}{2}\geq \frac{n-2+\sqrt{(n-2)^2+12n-8k-20}}{2},$$ 
with equality if and only if $n_{1}=2$ or $n_{1}=n-k-2$. If  $n_{1}=2$, then $G\cong \overrightarrow{G}_{n}^{k,2}$; if $n_{1}=n-k-2$, then $G\cong \overrightarrow{G}_{n}^{k,n-k-2}$, as required.

This completes the proof.
\begin{flushright}
$\square$
\end{flushright}
\textbf{Declaration of Completing Interest}\\

The authors declare that they have no known competing financial interests or personal relationships that could have appeared to influence the work reported in this paper.

\end{document}